
\input xy
\xyoption{all} 

\catcode`@=11  

\mag=1200

\hsize=140mm     \vsize=190mm
\hoffset=0mm     \voffset=0mm

\parindent=1cm

\catcode`\;=\active
\def;{\relax\ifhmode\ifdim\lastskip>\z@ 
\unskip\fi\kern.2em\fi\string;}

\catcode`\:=\active
\def:{\relax\ifhmode\ifdim\lastskip>\z@
\unskip\fi\penalty\@M\ \fi\string:}

\catcode`\!=\active
\def!{\relax\ifhmode\ifdim\lastskip>\z@
\unskip\fi\kern.2em\fi\string!}

\catcode`\?=\active
\def?{\relax\ifhmode\ifdim\lastskip>\z@
\unskip\fi\kern.2em\fi\string?}

\def\^#1{\if#1i{\accent"5E\i}\else{\accent"5E #1}\fi}
\def\"#1{\if#1i{\accent"7F\i}\else{\accent"7F #1}\fi}

\def\longlongrightarrow
{\relbar\joinrel\relbar\joinrel\relbar
\joinrel\relbar\joinrel\relbar\joinrel\rightarrow}

\def\CQFD{\quad\vbox{\hrule height 7pt
width 6pt}}

\frenchspacing
\catcode`\@=12 
\font\mega=cmbx10 scaled 2000


\nopagenumbers
\null\bigskip
\noindent{\mega Cohomologie des}
\medskip
\noindent{\mega alg\`ebres de Kr\oe necker
g\'en\'erales}
\bigskip\medskip

\noindent{\bf B. Bendiffalah}
(ben@math.univ-montp2.fr)

\noindent{\sevenrm Institut de
Math\'ematiques et de Mod\'elisation
de Montpellier (UMR CNRS 5149)}
\vskip -1mm
\noindent{\sevenrm
D\'epartement de Math\'ematiques,
Case 051}
\vskip -1mm
\noindent{\sevenrm
Universit\'e de Montpellier II,
Place Eug\`ene Bataillon}
\vskip -1mm

\noindent{\sevenrm MONTPELLIER 34095 Cedex 5,
FRANCE}
\bigskip

\noindent{\bf D. Guin}
(dguin@math.univ-montp2.fr)

\noindent{\sevenrm Institut de
Math\'ematiques et de Mod\'elisation
de Montpellier (UMR CNRS 5149)}
\vskip -1mm
\noindent{\sevenrm
D\'epartement de Math\'ematiques,
Case 051}
\vskip -1mm
\noindent{\sevenrm
Universit\'e de Montpellier II,
Place Eug\`ene Bataillon}
\vskip -1mm

\noindent{\sevenrm MONTPELLIER 34095 Cedex 5,
FRANCE}
\bigskip\medskip

{\bf Abstract for
``Cohomology of General Kr\"onecker Algebras''.}

The computation of the Hochschild
cohomology $H\!H^*T=H^*(T,T)$ of
a triangular algebra
$T=\big[{}^{A\,M}_{\,0\,\, B}\big]$
was performed in {\bf[BG2]}, by
the means of a certain triangular
complex. We use this result here
to show how $H\!H^*T$ splits in
little pieces whenever the bimodule
$M$ is decomposable. As an example,
we express the Hilbert-Poincar\'e serie
$\sum_{i=0}^\infty
dim_K(H\!H^iT_m)t^i$
of the ``general'' Kr\"onecker
algebra $T_m=\big[{}^{A\,M^m}_{\,0\ \,B}\big]$
as a function of $m\geq 1$
and those of $T$
(here the ground ring $K$ is a field
and $dim_KT<+\infty$).
The Lie algebra structure of
$H\!H^1T$ is also considered.
\bigskip

{\bf Sommaire.}
\smallskip

1. Introduction

\noindent 1.1 Cohomologie des
alg\`ebres triangulaires.

\noindent 1.2 R\'eductions.

\noindent 1.3 Alg\`ebres de Kr\oe necker
d'un bimodule.
\smallskip

2. Premier Th\'eor\`eme de R\'eduction

\noindent 2.1 Cohomologie triangulaire.

\noindent 2.2 Preuve du Lemme 1.

\noindent 2.3 Preuve du Th\'eor\`eme 1 (cf. 1.2.3).
\smallskip

3. Second Th\'eor\`eme de R\'eduction

\noindent 3.1 Cohomologie modifi\'ee.

\noindent 3.2 Suite exacte de Mayer-Vietoris.
\smallskip

4. L'Alg\`ebre de Lie $HH^1T$

\noindent 4.1  Preuve du Th\'eor\`eme 3.

\noindent 4.2 Scindage de $HH^1T$.
\smallskip

5. R\'ef\'erences Bibliographiques

\bigskip\vfill

{\bf Mots clefs :}
alg\`ebres triangulaires,
cohomologie de Hochschild,
suite exacte de Mayer-Vietoris,
alg\`ebres de Kr\oe necker g\'en\'eralis\'ees,
s\'erie de Hilbert-Poincar\'e.

{\bf Class. AMS :}
13D40,
16G60,
16E05,16E30,16E40,16E45,
20G05,20F40,
13J05.
\eject

\nopagenumbers
\null\vfill\eject
\pageno=3
\footline{\hfil
\bf---\ \oldstyle\folio\bf\ ---\hfil
}

\null\bigskip
\hfil\noindent{\bf COHOMOLOGIE DES ALG\`EBRES
DE KR\OE NECKER G\'EN\'ERALES}
\medskip
\hfil{\bf B. Bendiffalah} (Montpellier II) et
{\bf D. Guin} (Montpellier II)
\bigskip\bigskip

\goodbreak
\noindent\hfil
{\bf1. Introduction}
\medskip
De nombreux travaux ont montr\'e l'importance
du calcul de la cohomologie de \break Hochschild
des alg\`ebres associatives.
Il existe cependant peu de r\'esultats
g\'en\'eraux concernant les alg\`ebres
triangulaires ({\it cf. \bf \S1.1}).
Pourtant, la connaissance
de leur cohomologie est suffisante (par des
proc\'ed\'es de r\'ecurrence) pour obtenir
le calcul de la
cohomologie d'alg\`ebres tr\`es g\'en\'erales,
comme les alg\`ebres de 
poset (e.g. {\bf [M]}) ou les alg\`ebres
h\'er\'editaires arti\-niennes (e.g. {\bf [Har]}).
L'objectif de ce travail est d'\'etudier
l'incidence de la d\'ecomposabilit\'e
du bimodule $M$ sur la cohomologie de
Hochschild de l'alg\`ebre triangulaire.
En particulier, nous \'etudierons
la s\'erie de Hilbert-Poincar\'e
des alg\`ebres
de Kr\oe necker ``g\'en\'erales''
({\it cf. \bf \S1.3}).
\medskip

Un anneau commutatif $K$ est fix\'e
et $\otimes=\otimes_K$ ;
nous dirons ``module'' pour
tout module sur $K$.
\`A l'exception notable de
la section {\bf\S4}, o\`u il
est question d'alg\`ebres de Lie,
toutes nos
alg\`ebres sont des $K$-alg\`ebres
associatives unitaires.
Pour toute alg\`ebre $T$,
les ``$T$-modules''
sont des modules (unitaires)
\`a gauche sur $T$.
Pour les $T$-modules-\`a-droite,
nous utiliserons plut\^ot
l'alg\`ebre oppos\'ee $T^o$.

\medskip

\goodbreak
{\bf\S1.1 Cohomologie des alg\`ebres
triangulaires.}
\smallskip

Nous renvoyons \`a {\bf[CE]}
pour toutes les
questions basiques sur
la cohomologie de Hochschild.
Fixant deux alg\`ebres $A$ et $B$
et un $A\otimes B^o$-module
$M$ (un ``bimodule''),
on s'int\'eresse \`a l'alg\`ebre
``triangulaire'' ({\it e.g. \bf[ARS,III.\S2]}):
$$T=[{}^A_{\,0}{}^M_B]\ ,
\quad[{}^{a}_{0}{}^{m}_{\,b}]
[{}^{a'}_{\,0}{}^{m'}_{\,b'}]
=[{}^{aa'}_{\ 0}\,{}^{am'+mb'}_{\quad\ bb'}]
\ .\leqno\!\!(1.1.1)$$
Sa cohomologie,
$HH^*T:=H^*(T,T)$, n'est pas connue
dans le cas g\'en\'eral alors que
son homologie (de Hochschild) est ``triviale'' :
$H_*(T,T)=H_*(A,A)\oplus H_*(B,B)$
({\bf[Lo]}).

Hormis la section {\bf \S4},
nous supposerons syst\'ematiquement
que le $K$-module $T$ est projectif
(donc aussi $A$, $B$ et $M$) et nous
nous int\'eressons au
cas o\`u 
$M=\bigoplus_{i=1}^n M_i^{m_i}$ ($m_{i}\geq 1$),
sans aucune hypoth\`ese sur la
famille de
$A\otimes B^o$-modules $E=\{M_1,...,\!M_n\}$
(autre que les $M_i$ soient $K$-projectifs).
Nous montrons comment le 
calcul de
la cohomologie de l'alg\`ebre
$[{}^{A\,M}_{\,0\,\,B}]$ se d\'eduit
de ceux des alg\`ebres
$[{}^{A\,M_i}_{\,0\,\,B}]$,
\`a l'aide de deux sortes de
r\'eductions ({\it cf. \bf\S1.2}).
La premi\`ere (Th\'eor\`eme 1)
ram\`ene
le calcul de $HH^*[{}^A_{\,0}{}^{M}_B]$
\`a celui de $HH^*[{}^A_{\,0}{}^{M'}_B]$,
o\`u $M'=\bigoplus_iM_i$ : nous en
d\'eduisons la cohomologie des alg\`ebres
de Kr\oe necker g\'en\'erales
({\it cf. \bf\S1.3}).
La seconde r\'eduction
(Th\'eor\`eme 2) est
une famille de suites exactes longues
de type Mayer-Vietoris ramenant le calcul
de $HH^*[{}^A_{\,0}{}^{M'}_B]$ \`a
ceux des $HH^*[{}^{A\,M_i}_{\,0\ B}]$,
$1\leq i\leq n$.

Ce travail fait suite \`a
{\bf [BG2]}, o\`u la cohomologie
de $T$ est exprim\'ee
\`a l'aide d'un
certain complexe ``triangulaire''
($M$ est suppos\'e $K$-projectif).
Ce r\'esultat, que nous rappelons
dans {\bf\S2.1} pour le confort du lecteur,
avait permis entre autres choses
de retrouver la g\'en\'eralisation
de {\bf[Ci]} et {\bf [MP]}
de la suite exacte longue
de D. Happel {\bf[Hap]} 
({\it cf.} aussi
{\bf[GMS]},{\bf[GS]},{\bf[GG]},{\bf[K]}):
$$\cdots\longrightarrow
Ext^{*-1}_{A\otimes B^o}(M,M)
\buildrel{l^*_M}\over\longrightarrow
HH^*[{}^A_{\,0}{}^{M}_B]
\buildrel{r^*_M}\over\longrightarrow
\matrix{HH^*A\cr HH^*B\cr}
\buildrel{d^*_M}
\over\longrightarrow
Ext^*_{A\otimes B^o}(M,M)
\longrightarrow\cdots\leqno\!\!(1.1.2)$$
o\`u $r^*_M$ est un morphisme
d'alg\`ebres gradu\'ees
(la superposition signifie un
produit direct).
Incidemment ({\it cf.} 1.2.15, avec $N=0$),
nous retrouvons que $r^1_M$ est un morphisme
d'alg\`ebres de Lie (sans hypoth\`ese sur $T$).
\bigskip

{\bf\S1.2 R\'eductions.}

Le r\'esultat suivant,
d\'emontr\'e dans {\bf\S2.2}, est une
g\'en\'eralisation de 1.1.2
($l^*_{0,M}=l^*_M$ et $r^*_{0,M}=r^*_M$).

\proclaim 1.2.1 Lemme 1.
Pour tout $A\otimes B^o$-module
$N$ ($K$-projectif), on a une
suite exacte longue :
$$\cdot\cdot\!\longrightarrow\!
HH^*\big[{}^A_{\,0}
\,{}^{N}_{\,B}\big]
\buildrel{d^{*-1}_{N\!,M}}
\over\longrightarrow\!
\matrix{Ext^{*-1}_{A\otimes B^o}(N\!,M)\cr
Ext^{*-1}_{A\otimes B^o}(M,M)\cr
Ext^{*-1}_{A\otimes B^o}(M,N)\cr}
\buildrel{l^*_{N\!,M}}
\over\longrightarrow
HH^*\big[{}^A_{\,0}
\,{}^{M\oplus N}_{\ \ \,B}\big]
\buildrel{r^*_{N,M}}
\over\longrightarrow
HH^*\big[{}^A_{\,0}
\,{}^{N}_{B}\big]
\!\!\longrightarrow\!
\cdot\cdot\leqno \!\!(1.2.2)$$
o\`u $l^*_{N,M}$ est induit par
$l^*_{M\oplus N}$ ({\it cf.} 1.1.2)
et $r^*_{L,N}\circ
r^*_{N\oplus L,M}=
r^*_{L,N\oplus M}$.

\noindent
Ce lemme permet un premier d\'evissage 
de $H\!H^*[{}^A_{\,0}{}^{M}_B]$ si
$M=\bigoplus_{i=1}^nM_i^{m_i}$, pour une
famille
de $A\!\otimes\! B^o$-modules $K$-projectifs
 $E\!=\!\{M_1,..,M_n\}$ et $m_i\geq 1$.
Posant $M'\!=\!\bigoplus_{i=1}^n M_i$,
le th\'eor\`eme suivant
(d\'emontr\'e dans {\bf\S2.3})
nous autorise \`a faire abstraction
des multiplicit\'es $m_i$.

\proclaim 1.2.3 Th\'eor\`eme 1.
Nous avons une suite exacte
courte gradu\'ee scind\'ee de modules :
$$0\longrightarrow
\prod_{i,j}Ext^{*-1}_{A\otimes B^o}
(M_i,M_j)^{m_im_j-1}
\buildrel{l^*}\over\longrightarrow
HH^*[{}^A_{\,0}{}^{M}_B]
\buildrel{r^*}\over
\longrightarrow HH^*[{}^A_{\,0}{}^{M'}_B]
\longrightarrow 0\ ,\leqno \!\!(1.2.4)$$
o\`u les morphismes $l^*$ et $r^*$
sont d\'efinis par le Lemme 1.

Pour compl\'eter notre programme,
il nous reste
\`a exprimer
$HH^*\big[{}^A_{\,0}{}^{M'}_B\big]$
en fonction des $HH^*\big[{}^A_{\,0}
{}^{M\!\scriptscriptstyle i}_{\,B}\big]$.
C'est un objectif que nous
r\'ealisons au moyen de suites
exactes longues
qui s'apparentent \`a celle de
Mayer-Vietoris (Th\'eor\`eme 2).
Pour les d\'ecrire,
il sera utile, au regard du
r\'esultat suivant, de modifier
la cohomologie employ\'ee. 

\proclaim 1.2.5 Lemme 2.
Le module
$\prod_{i\not=j}Ext^{*-1}_{A\otimes B^o}
(M_i,M_j)$
est (isomorphe par $l^*_{M'}$ \`a) un
sous-module de
$HH^*\big[{}^A_{\,0}{}^{M'}_B\big]$.

En fait ({\it cf. \bf\S3.1}), le sous-module
de $HH^*\big[{}^A_{\,0}{}^{M'}_B\big]$
pr\'ecis\'e au Lemme 2 est en facteur
direct et poss\`ede un suppl\'ementaire
particulier, dont on donnera une
d\'efinition intrins\`eque
avec 3.1.5. 

\proclaim 1.2.6 D\'efinition 1.
Nous notons ${\cal H}^*E$
ce suppl\'ementaire canonique
(de l'image par $l^*_{M'}$) de $\prod_{i\not=j}
Ext^{*-1}_{A\otimes B^o}
(M_i,M_j)$ dans
$HH^*\big[{}^A_{\,0}{}^{M'}_B\big]$.
Par convention,
${\cal H}^*(\phi)=
HH^*\big[{}^A_{\,0}{}^{\,0}_B\big]$.

\noindent
Dans la suite, nous aurons \`a consid\'erer des
sous-familles $F$ de la famille $E$  et nous
noterons $\bar F$ la somme directe de ses
\'el\'ements ;  on pose $\bar \phi = 0$.
Avec cette notation, on a
$M' = \bar E$.
Ceci permet de formuler
imm\'ediatement la g\'en\'eralisation
suivante de 1.2.2  ({\it cf. \bf\S3.1}).
 
\proclaim 1.2.7 Lemme 3.
Pour toute partie $F\subset E$,
nous avons une suite exacte longue :
$$\cdots\longrightarrow
\!\!\!\!\prod_{L\in E\backslash F}\!\!\!
Ext^{*-1}_{A\otimes B^o}(L,L)
\buildrel{\ell^*_{F,E}}\over\longrightarrow
{\cal H}^*E
\buildrel{\rho^*_{F,E}}\over
\longrightarrow {\cal H}^*F
\buildrel{\delta^*_{F,E}}\over\longrightarrow
\!\!\!\!\prod_{L\in E\backslash F}\!\!\!
Ext^{*}_{A\otimes B^o}(L,L)
\longrightarrow\cdots\leqno\!\!(1.2.8)$$
o\`u $\ell^*_{F,E}$, $\rho^*_{F,E}$
et $\delta^*_{F,E}$ sont induits par
$l^*_{\bar F,\bar E}$,
$r^*_{\bar F,\bar E}$ et
$d^*_{\bar F,\bar E}$
({\it cf.} 1.2.2).

\noindent En particulier, pour 
$G\subset F\subset E$,
nous avons 
$\rho^*_{G,F}\circ \rho^*_{F,E}=\rho^*_{G,E}$
(fonctorialit\'e).

On retrouve 1.2.2 avec $E=\{M,N\}$ et
$F=\{N\}$ ; 
de m\^ eme 1.1.2 s'obtient avec
$E=\{M\}$ et $F=\phi$.
Le th\'eor\`eme suivant d\'ecoule
du Th\'eor\`eme 4, {\bf\S3.2}.

\proclaim 1.2.9 Th\'eor\`eme 2.
Pour toutes
parties $F,G\subset E$, nous avons
une suite exacte longue
$$\cdot\cdot\!\longrightarrow\!
{\cal H}^*(F\cup G)
\buildrel{\textstyle
{}^{\rho^*_{F,F\cup G}}_{\rho^*_{G,F\cup G}}}
\over\longlongrightarrow
\matrix{{\cal H}^*F\cr {\cal H}^*G\cr}
\buildrel{\rho^*_{F\cap G,F}-\rho^*_{F\cap G,G}}
\over\longlongrightarrow {\cal H}^*(F\cap G)
\longrightarrow{\cal H}^{*+1}(F\cup G)
\!\longrightarrow\!\cdot\cdot
\leqno \!\!(1.2.10)$$

Dor\'enavant, le calcul de
$HH^*\big[{}^A_{\,0}{}^{M}_B\big]$
se ram\`ene, dans les bons cas de
d\'ecomposabilit\'e de $M$ 
({\it e.g.} longueur finie),
\`a ceux des
$A\otimes B^o$-modules
ind\'ecomposables.
Pour \^etre vraiment complet, il
faudrait aussi pouvoir suivre toute
la stucture alg\'ebrique de
$HH^*\big[{}^A_{\,0}{}^{M}_B\big]$
dans nos r\'eductions.
\`A ce propos, nous savons
que ${\cal H}^*E$ est muni
d'une structure d'alg\`ebre
(isomorphe \`a une alg\`ebre quotient
de $HH^*[{}^{A\,\bar E}_{\,0\, B}]$, 
{\it cf.} {\bf[GMS]} et {\bf[GS]}).
De l\`a, notre premi\`ere question :
\smallskip

\noindent$\!\!(1.2.11)$ {\sl Les morphismes
$\rho^*_{F,E}$
sont-ils des morphismes d'alg\`ebres ?}
\smallskip

\noindent
Question que nous laissons sans r\'eponse ;
cependant, nous ne pensons pas que cela soit
vrai pour le morphisme $r^*_{N,M}$
de 1.2.2 ;
le cas  $r^*_M=r^*_{0,M}$ de 1.1.2
serait donc tr\`es exceptionnel.

Dans le cas d'une alg\`ebre quelconque $T$,
la structure d'alg\`ebre de Lie de $HH^1T$
est extr\^emement difficile \`a \'elucider ;
si $T$ est une alg\`ebre monomiale
de dimension finie ($K$ est un corps
alg\'ebriquement clos), l'on dispose
d'une tr\`es jolie interpr\'etation
g\'eom\'etrique associ\'ee au groupe de Lie
$Aut(T)$ (les automorphismes de $T$)
et la situation s'en trouve mieux comprise
({\it e.g. \bf[GAS]}, {\bf[S]}). Dans le cadre
des alg\`ebres triangulaires, il se trouve
que le module ${\cal H}^1E$ d\'efini en 1.2.6
est une sous-alg\`ebre de Lie de $\!H\!H^1\big[
{}^{A\,\bar E}_{\,0\, B}\big]$. De l\`a
une seconde question, aussi naturelle,
\`a laquelle
nous r\'epondons positivement :
\smallskip

\noindent$\!\!(1.2.12)$ {\sl Les morphismes
$\rho^1_{F,E}$ de 1.2.8
sont-ils des morphismes d'alg\`ebres de Lie ?}
\smallskip

\noindent
(Nous ne pensons pas que cela soit
g\'en\'eralement vrai pour $r^1_{N,M}$.)
Par fonctorialit\'e de $\rho^1$, il suffit
de r\'epondre \`a 1.2.12 dans le cas
particulier de l'alg\`ebre de Lie 
$HH^1\big[
{}^{A\,M\oplus N}_{\,0\quad B}\big]$ :

\proclaim 1.2.13 Th\'eor\`eme 3.
Il existe deux sous-alg\`ebres de Lie,
${\cal H}^1(M,N)$ et ${\cal H}^1(N,M)$,
v\'erifiant
$${\cal H}^1(M,N)\!+\!{\cal H}^1(N,M)
\!=\!H\!H^1\big[
{}^{A\,M\oplus N}_{\,0\quad B}\big]
\quad et\quad
{\cal H}^1(M,N)\cap {\cal H}^1(N,M)
\!=\!{\cal H}^1\{M,N\}\ .
\leqno\!\!(1.2.14)$$
Par restriction,
$r^1_{M,N}$ induit des morphismes
d'alg\`ebres de Lie de m\^eme image
que $r^1_{M,N}$ :
$$\eqalign{\rho^1_{M,(M,N)}&:{\cal H}^1(M,N)
\longrightarrow HH^1\big[
{}^{A\,M}_{\,0\ B}\big]\ ,\cr
\rho^1_{M,(N,M)}&:{\cal H}^1(N,M)
\longrightarrow HH^1\big[
{}^{A\,M}_{\,0\ B}\big]\ ,\cr
\rho^1_{M,\{M,N\}}&:{\cal H}^1\{M,N\}
\longrightarrow HH^1\big[
{}^{A\,M}_{\,0\ B}\big]\ .\cr}
\leqno\!\!(1.2.15)$$

\noindent
Dans la section {\bf \S4}, nous donnons
une condition n\'ecessaire et suffisante
(Th\'eor\`eme 5), portant sur $N$, pour
avoir la surjectivit\'e du morphisme
$r^1_{M,N}$, donc des morphismes
d'alg\`ebres de Lie 1.2.15.
En particulier, nous en d\'eduirons
une suite exacte :
$$0\!\longrightarrow\!
Z\!\big[{}^{A\,M\oplus N}_{\,0\quad B}\big]
\!\buildrel{r^0_{M,N}}\over\longrightarrow\!
Z\!\big[{}^{A\,M}_{\,0\, B}\big]
\!\longrightarrow\! End_{A\otimes B^o}N
\!\longrightarrow\!
{\cal H}^1\!\{M,\!N\}
\!\buildrel{\rho^1_{M,\{\!M,\!N\}}}\over
{-\!-\!\!\!\longrightarrow}\!
H\!H^1\!\big[{}^{A\,M}_{\,0\, B}\big]
\!\longrightarrow\! 0.\leqno\!\!(1.2.16)$$
\smallskip

Dans le paragraphe suivant nous
appliquons nos r\'esultats \`a une
g\'en\'eralisation des alg\`ebres
de Kr\oe necker.
\medskip

\goodbreak
{\bf\S1.3 Alg\`ebres de
Kr\oe necker d'un bimodule.}

L'alg\`ebre de Kr\oe necker ${\bf K}_2$ est
l'alg\`ebre du carquois
$\bullet
{}^{\displaystyle\longrightarrow
}_{\displaystyle\longrightarrow}
\bullet$ (deux sommets et deux fl\`eches)
et nous avons un isomorphisme
d'alg\`ebres ${\bf K}_2\cong
\big[{}^K_{\,0}{}^{K^2}_{\,K}\big]$.
Si $K$ est un corps, c'est l'exemple
le plus simple (et le plus \'etudi\'e)
d'alg\`ebre dont le type de repr\'esentation
n'est pas fini et qui ne soit pas sauvage
({\it e.g. \bf[ARS,VIII.\S7]}) ;
pour $m\geq 3$, l'alg\`ebre de
Kr\oe necker ``g\'en\'eralis\'ee''
${\bf K}_m=
\big[{}^K_{\,0}{}^{K^m}_{\,K}\big]$
est sauvage.

Pour $m\geq 1$, on sait que le $K$-module 
$H\!H^*\big[{}^K_{\,0}{}^{K^m}_{\,K}\big]$
est libre et l'on a
$$rang_K\big(H\!H^i\big[{}^{K\, K^m}_{\,0\ K}
\big]\big)
=\cases{\quad\ 1& si $i=0$ ;\cr
m^2-1& si $i=1$ ;\cr
\quad\  0&sinon.\cr}\leqno \!\!(1.3.1)
$$
({\it e.g. \bf[BG1]}, l'\'etude
y est faite pour un cardinal $m$
quelconque, m\^eme infini, en connection
avec les cat\'egories ``muscl\'ees'' ;
{\it cf.} aussi 4.2.10).
Nous d\'esirons
g\'en\'eraliser la formule 1.3.1
en associant \`a tout $A\otimes B^o$-module
$M$ (avec $A$, $M$ et $B$, $K$-projectifs),
des ``alg\`ebre de Kr\oe necker g\'en\'erales''
$[{}^A_{\,0}\,{}^{M^m}_{\,B}\big]$
(alg\`ebres triangulaires), $m\geq 1$.
D'apr\`es le Th\'eor\`eme 1,
nous avons une suite exacte courte
gradu\'ee scind\'ee :
$$0\longrightarrow
Ext^{*-1}_{A\otimes B^o}(M,M)^{m^2-1}
\longrightarrow
HH^*\big[{}^{A\,M^m}_{\,0\ B}\big]
\longrightarrow
HH^*[{}^A_{\,0}{}^M_B]
\longrightarrow 0\ .\leqno \!\!(1.3.2)$$

Pour $M\not=0$ et $m\geq 2$,
ce r\'esultat prouve, par exemple, que nous
avons toujours
$HH^1\big[{}^{A\,M^m}_{\,0\ B}\big]\not=0$.
Si $K$ est un corps et $A$, $B$ et 
$N$ sont des $K$-espaces vectoriels
de dimensions finies, nous avons la
s\'erie de 
Hilbert-Poincar\'e\footnote{${}^{(\dagger)}$}
{$\chi_A=\sum_{i=0}^\infty
t^i\dim_K(H\!H^iA)\in {\bf Z}[[\,t\,]]$} :
$$\chi_{\left[{}^A_{\,0}
{}^{M^{\!m}}_{\,B}\right]}
=\chi_{\left[{}^A_{\,0}{}^{M}_{B}\right]}
+t(m^2-1)\Xi\ ,\leqno \!\!(1.3.3)$$
o\`u $\Xi=\sum_{i=0}^\infty
t^i\dim_KExt^i_{\!A\otimes B^o}(M,M)$. 
Si $car (K)=p>0$, nous en d\'eduisons
la $p$-p\'eriodicit\'e de l'application 
$m\longmapsto \chi'(m)$,
o\`u $\chi'(m)\!\in\! {\displaystyle {\bf Z}/
{\displaystyle p\bf Z}}[[t]]
\!\subset\! K[[t]]$
est la classe  de
$\chi_{\left[{}^A_{\,0}{}^{M^{\!m}}_{\,B}\right]}$
modulo $p$ ; $\chi'(m)$ ne d\'epend que
de la classe de $m^2$ modulo $p$
et en particulier :
$$\chi'(p-1)=
\chi'(1)=
\chi'(p+1)\ .
\leqno \!\!(1.3.4)$$

Le corollaire suivant g\'en\'eralise
un r\'esultat de
{\bf[BG2]}, rappel\'e dans {\bf\S2.1}
({\it cf.} 2.1.11).

\proclaim 1.3.5 Corollaire 1.
Soit une alg\`ebre $A$ dont le module
sous-jacent est projectif.
Pour tout $A$-module projectif $M$,
de type  fini (sauf si $K$ est
semi-simple) et pour tout entier
$m\geq 1$, nous avons une suite exacte
courte scind\'ee gradu\'ee :
$$0\longrightarrow
(HH^{*-1}B)^{m^2-1}
\longrightarrow
HH^*[{}^A_{\,0}{}^{M^m}_{\,B}]
\buildrel{r^*}\over
\longrightarrow HH^*A
\longrightarrow 0\ ,\leqno \!\!(1.3.6)$$
o\`u $B=(End_AM)^o$
(ainsi $M$ est un $A\otimes B^o$-module)
et $r^*$ est un morphisme
d'alg\`ebres.

\noindent 
En particulier (m\^emes
hypoth\`eses) : si $K$ est un corps
et si les $K$-espaces vectoriels
$A$ et $N$ sont de dimension finie,
nous avons 
$\chi_{\left[{}^A_{\,0}{}^{M^{\!m}}_{\,B}\right]}
=\chi_A+t(m^2-1)\chi_B$ et la s\'erie 
$\chi_{\left[{}^A_{\,0}{}^{M^{\!m}}_{\,B}\right]}$
est polynomiale (resp. rationnelle)
si $\chi_A$ et $\chi_B$ le sont.
\bigskip

\noindent\hfil{\bf2. Premier Th\'eor\`eme
de R\'eduction}
\smallskip

Dans le paragraphe {\bf \S2.2},
nous explicitons des sous-modules
canoniques de la cohomologie
d'une alg\`ebre triangulaire
dont le bimodule
est d\'ecomposable ; la suite
exacte longue 1.2.2
est d\'emontr\'ee.
On obtient ensuite le Th\'eor\`eme 1
de {\bf\S1.2} en scindant cette
suite exacte longue, gr\^ace au
crit\`ere d\'efini par la Proposition 4
(ou 4 {\it bis}) de {\bf \S2.3}.

Commen\c cons par rappeler
la d\'efinition de la cohomologie
de Hochschild et la construction du complexe
triangulaire introduit dans {\bf [BG2]}.
\smallskip

\goodbreak
{\bf\S2.1 Cohomologie triangulaire.}

Pour toute alg\`ebre $T$, nous avons
un ``bar-complexe'' ({\it cf. \bf[Lo]})
$C^T_*=\bigoplus_{n\geq 0}C^T_n$,
avec $C^T_n=T\otimes\cdots\otimes T$
($n+2$ facteurs $T$) et de
diff\'erentielle
$$a_0\otimes\cdots\otimes a_{n+1}
\longmapsto
\sum_{i=0}^n(-1)^ia_0\otimes\cdots
\otimes a_ia_{i+1}\otimes\cdots\otimes a_{n+1}\ .
\leqno\!\!(2.1.1)$$
C'est un complexe de $T\otimes T^o$-modules
que nous pouvons augmenter en un complexe
acyclique $C^{T+}_*$ avec la
multiplication $T\otimes T\longrightarrow T$,
$C_{-\!1}^{T+}=T$.
Le bar-complexe permet de d\'efinir, pour tout
$T\otimes T^o$-module $\Lambda$,
le complexe (croissant) de Hochschild
$C^*(T,\Lambda)=Hom_{T\otimes T^o}(C^T_*,\Lambda)$
et, donc, la cohomologie de
Hochschild $H^*(T,\Lambda)$.
La cohomologie de $C^*(T)=C^*(T,T)$,
not\'ee $HH^*T=H^*(T,T)$,
est une alg\`ebre associative gradu\'ee
et $HH^1T$ est naturellement munie
une structure d'alg\`ebre de Lie dont
nous reparlerons dans {\bf \S4.1}.
Si $T$ est un $K$-module projectif,
il existe un isomorphisme
$H^*(T,\Lambda)\cong
Ext^*_{T\otimes T^o}(T,\Lambda)$ qui,
pour $\Lambda=T$, est un isomorphisme
d'alg\`ebres gradu\'ees (produit de Yoneda).

Soient \`a pr\'esent deux alg\`ebres
$A$ et $B$ et un $A\!\otimes\! B^o$-module $M$,
tous $K$-projectifs.
Le quotient de $(C_*^{A+}\!\otimes\! C_*^{B^o+})[-1]$
(suspension du produit des complexes) par son
sous-complexe ponctuel $\!A\otimes\! B[-1]$
(concentr\'e en degr\'e $-1$) est un
complexe positif de $A\!\otimes\! A^o\!\otimes\!
B\!\otimes\! B^o$-modules,
not\'e $C^{A,B^o}_*\!\!\!\!,\ $
et nous d\'efinissons
$C_*^{A,B^o}\!\!M=
C^{A,B^o}_*\!\!\otimes_{A\otimes B^o} M$.
En particulier : $C^{A}_*M=
C^{A}_*\otimes_{A} M$ et 
$C^{B^o}_*\!M=M\!\otimes_{B}\! C^{B}_*$
sont des sous-complexes
disjoints de $C_*^{A,B^o}\!\!M$ :
nous notons
$i^{A}_{M,N}$ et $i^{B^o}_{M,N}$
ces injections.
Les morphismes canoniques
d'alg\`ebres induits par la structure
de $A\otimes B^o$-module de $M$,
$\alpha_M:A\longrightarrow End_{B^o}M$
et $\beta_M:B\longrightarrow (End_{A}M)^o$,
induisent des
morphismes de complexes
$\alpha^*_M:
C^*(A)\longrightarrow C^*(A,End_{B^o}M)$
et
$\beta^*_M:
C^*(B)\longrightarrow C^*(B,End_{A}M)$.

\proclaim 2.1.2 D\'efinition 2.
Pour tout $A\otimes B^o$-module $N$,
le complexe ``triangulaire de $M$
\`a coefficient dans $N$'' est
$C^*_{tri}(M,N)=Hom_{A\otimes B^o}
(C^{A,B^o}_*\!M,N)$, not\'e
$C^*_{tri}(M)$ si $M=N$.

Avec les isomorphismes d'adjonction
$C^*(A,Hom_{B^o}(M,N))\cong
Hom_{A\otimes B^o}(C^A_*M,N)$
et
$C^*(B,Hom_{A}(M,N))\cong
Hom_{A\otimes B^o}(C^{B^o}_*M,N)$,
nous d\'eduisons un morphisme
surjectif :
\smallskip

\noindent $\!\!$(2.1.3)\ \
$i^*_{M,N}\!=\!(i^{A*}_{M,N},i^{B^o\!\!*}_{M,N}):
C^*_{tri}(M,N)\longrightarrow
C^*(A,Hom_{B^o}(M,N))\times
C^*(B,Hom_{A}(M,N))$.
\medskip

\noindent{\bf 2.1.4 Remarque.}
Le complexe $\Sigma_{M,N}^*$,
c\^one du morphisme $i^*_{M,N}$,
est quasi-isomorphe au complexe
$Ker\,i^*_{M,N}$, dont
l'homologie
est isomorphe au module gradu\'e
$Ext^{*-1}_{A\otimes B^o}(M,N)$.

\proclaim 2.1.5 D\'efinition 3.
Pour $T\!=\!\big[{}^A_{\,0}{}^M_B\big]$,
on note  ${\cal C}^*\{M\}$ le
c\^one du morphisme de
complexes\footnote{${}^{(\dagger)}$}
{\rm ${\cal C}^*\{M\}$ est not\'e
$\tilde C^*(T,T)$ dans {\bf[BG2]}}
$$\lambda^*_M:
C^*(A)\times C^*_{tri}(M)\times C^*(B)
\longrightarrow
C^*(A,End_{B^o}M)\times C^*(B,End_{A}M),
\leqno\!\!(2.1.6)$$
dont les composantes
\'eventuellement non nulles
sont $\alpha^*_M$, $i^{A*}_M=i^{A*}_{M,M}$,
 $i^{B^o\!\!*}_M=i^{B^o\!\!*}_{M,M}$
et $\beta^*_M$.

Nous en d\'eduisons une suite exacte
courte de complexes
($\Sigma^*_M=\Sigma^*_{M,M}$) :
$$0\longrightarrow \Sigma^*_M
\buildrel{l_M}\over\longrightarrow
 {\cal C}^*\{M\}
\buildrel{r_M}\over\longrightarrow
C^*(A)\times  C^*(B)
\longrightarrow 0\ .\leqno\!\!(2.1.7)$$

\proclaim 2.1.8 Th\'eor\`eme [BG2, 3.4.4].
Nous avons
une \'equivalence homotopique de
complexes :
$$C^*(T)\sim {\cal C}^*\{M\}
\leqno\!\!(2.1.9)$$
et la suite exacte longue
1.1.2 provient de la suite exacte
courte de complexes 2.1.7.

\noindent Exemple d'application :
si $M$ est un $A$-module
projectif, nous avons une
suite exacte longue
$$\cdots\longrightarrow HH^{*-1}A
\longrightarrow
H(c\hat one\,\beta^*_M)\longrightarrow
HH^*[{}^A_{\,0}{}^M_{B}]
\buildrel{r^*}\over\longrightarrow HH^*A
\longrightarrow\cdots\leqno\!\!(2.1.10)$$
o\`u $r^*$ est un morphisme
d'alg\`ebres gradu\'ees et
$H(c\hat one\,\beta^*_M)$
est l'homologie du c\^one du morphisme
$C^*(B,B)\longrightarrow C^*(B,End_AM)$
que $\beta_M$ induit. En particulier,
nous avons :
\proclaim 2.1.11 Corollaire [BG2, 1.3.7].
Si $A$ est une alg\`ebre
$K$-module projectif et
$M$ est un $A$-module
projectif de type fini,
nous avons des isomorphismes
d'alg\`ebres gradu\'ees :
$$HH^*[{}^A_{\,0}{}^M_{B}]\cong 
HH^*A\quad\hbox{\sl 
(induit par $r^*_M$) et}\quad
HH^*B\cong Ext^{*}_{A\otimes B^o}(M,M)
\ ,\leqno\!\!(2.1.12)$$
o\`u $B=(End_AM)^o$
(ainsi $M$ est bien un $A\otimes B^o$-module).

\noindent En fait l'hypoth\`ese de
finitude sur $M$ n'est l\`a que pour
assurer que l'alg\`ebre $B$ est un module
$K$-projectif :
elle est redondante si, par exemple,
l'anneau $K$ est semi-simple.
\medskip

Clairement, le Corollaire 1 de {\bf\S1.3}
s'obtient en mettant ensemble
la suite scind\'ee 1.3.2 
avec les isomorphismes 2.1.12.
La suite exacte courte gradu\'ee 1.3.6
provient de la suite exacte longue 2.1.10 :
$H(c\hat one\,\beta^*_{M^m})
\cong H^{*-1}(B,End_AM)^{m^2-1}\cong
(HH^{*-1}B)^{m^2-1}$.
\medskip

{\bf\S2.2 Preuve du Lemme 1 (\S1.2).}

Soient deux alg\`ebres  $A$ et $B$
et soient deux $A\otimes B^o$-modules
$M$ et $N$ (tous suppos\'es $K$-projectifs).
Le r\'esultat suivant est un cas
particulier du Lemme 2 de {\bf\S1.2}
(prendre $n=2$).
\proclaim 2.2.1 Proposition 1.
Le morphisme $l^*_{M\oplus N}$
({\it cf.} 1.1.2) est injectif
sur le sous-module
$$Ext^{*-1}_{A\otimes B^o}(M,N)\times
Ext^{*-1}_{A\otimes B^o}(N,M)
\subset Ext^{*-1}_{A\otimes B^o}
(M\oplus N,M\oplus N)\leqno\!\!(2.2.2)$$
et son image 
est un facteur direct de
$HH^*[{}^A_{\,0}\,{}^{M\oplus N}_{\ \ B}]$.

\noindent
Autrement dit : il existe
un sous-module gradu\'e
${\cal H}^*\{M,N\}\subset
HH^*\big[{}^A_{\,0}\,{}^{M\oplus N}_{\ \ B}\big]$
(d\'efini dans la preuve ci-dessous)
et un isomorphisme
$$HH^*[{}^A_{\,0}\,{}^{M\oplus N}_{\ \ B}]
\cong {\cal H}^*\{M,N\}\times
Ext^{*-1}_{A\otimes B^o}(M,N)\times
Ext^{*-1}_{A\otimes B^o}(N,M)\ .
\leqno\!\!(2.2.3)$$
\medskip

\noindent{\bf 2.2.4 Preuve.}
Consid\'erons le complexe 
${\cal C}^*\{M\oplus N\}$
(c\^one du morphisme
$\lambda_{M\oplus N}^*$).
Dans la d\'ecomposition de Peirce
$$End_{A}(M\!\oplus\! N)\!\cong\!
\left[\matrix{\!\!End_AM
&\!\!\!\!\!\!\!H\!om_A(N,M)\!\!\cr
\!\!H\!om_A(M,N)
&\!\!\!\!\!\!\!End_AN\!\!\cr}\right]\ ,$$
ainsi que celle de l'alg\`ebre
$End_{B^o}(M\!\oplus\! N)$,
il est essentiel de noter que
l'image de $\alpha_{M\oplus N}^*$
(resp. $\beta_{M\oplus N}^*$)
est contenue dans le sous-complexe
$$C^*(A,End_{B^o}M)\times C^*(A,End_{B^o}N)
\quad({\rm resp.}\ C^*(B,End_{A^o}M)
\times C^*(B,End_{A^o}N)\ )\leqno\!\!(2.2.5)$$
du complexe
$C^*\big(A,End_{B^o}(M\oplus N)\big)$
(resp. $C^*\big(B,End_{A^o}(M\oplus N)\big)$).
Nous en d\'eduisons une
d\'ecomposition en somme directe
de complexes :
$${\cal C}^*\{M\oplus N\}
={\cal C}^*\{M,N\}
\oplus \Sigma^*\{M,N\}\ ,
\leqno\!\!(2.2.6)$$
o\`u $\Sigma^*\{M,N\}\!=\!\Sigma^*_{M,N}\!\times\!
\Sigma^*_{N,M}$ et ${\cal C}^*\{M,N\}$
est le complexe c\^one du morphisme
de complexes :
$$\lambda^*_{M,N}:
\matrix{{}^{\quad C^*(A)}_{C^*_{tri}(M,M)}\cr
{}^{C^*_{tri}(N,N)}_{\quad C^*(B)}\cr}
\longlongrightarrow \matrix{
{}^{C^*(A,End_{B^o}M)}_{C^*(A,End_{B^o}N)}\cr
{}^{C^*(B,End_{A}M)}_{C^*(B,End_AN)}\cr}
\leqno\!\!(2.2.7)$$
(la superposition encode
une somme directe de complexes)
ayant pour seules composantes \'eventuellement
non nulles : $\alpha_M^*$, $\alpha_N^*$,
$\beta_M^*$, $\beta_N^*$, $i^{A*}_{M}$,
$i^{B^o*}_{M}$,  $i^{A*}_{N}$ et
$i^{B^o*}_{N}$.
Notant ${\cal H}^*\{M,N\}$ la
cohomologie du complexe
${\cal C}^*\{M,N\}$, nous obtenons
la Proposition 1 avec 2.1.4.
\CQFD\medskip

La Proposition 2 suivante, alli\'ee
\`a la Proposition 1, prouve le
Lemme 1 de {\bf\S1.2}.

\proclaim 2.2.8 Proposition 2. Nous avons
une suite exacte longue :
$$\cdot\cdot\!\longrightarrow
Ext^{*-1}_{A\otimes B^o}(M,M)
\!\buildrel{\ell^*_{N,M}}
\over\longrightarrow\!
{\cal H}^*\{M,N\}
\!\buildrel{\rho^*_{N,M}}
\over\longrightarrow\!
HH^*\big[{}^A_{\,0}
\,{}^{N}_{B}\big]
\!\buildrel{\delta^*_{N,M}}
\over\longrightarrow\!
Ext^{*}_{A\otimes B^o}(M,M)
\!\longrightarrow\!
\cdot\cdot\leqno\!\!(2.2.9)$$
o\`u $\ell^*_{N,M}$ est induit par
$l^*_{M\oplus N}$ et $\rho^*_{N,M}$
(induit par $r^*_{M\oplus N}$) v\'erifie :
$r^*_N\circ \rho^*_{N,M}=
r^*_M\circ \rho^*_{M,N}$.

\noindent{\bf 2.2.10 Preuve.}
Reprenons les notations de 2.2.4.
Nous avons un morphisme
de complexes
$\rho_{N,M}:
{\cal C}^*\{M,N\}\longrightarrow
{\cal C}^*\{N\}$,
d\'efini de mani\`ere \'evidente
en envoyant toutes les composantes
d\'ecrites en 2.2.7 sur elles-m\^emes,
except\'ees
$C^*_{tri}(M)$, $C^*(A,End_{B^o}M)$
et $C^*(B,End_{A}M)$ sur lesquelles
$\rho_{N,M}$ est nulle.
Il s'agit d'un morphisme 
surjectif, de noyau $\Sigma^*_M$.
De la suite exacte courte de complexes
$$0\longrightarrow\Sigma^*_{M}
\buildrel{\ell_{N,M}}\over
\longrightarrow {\cal C}^*\{M,N\}
\buildrel{\rho_{N,M}}\over
\longrightarrow {\cal C}^*\{N\}
\longrightarrow 0\ ,
\leqno\!\!(2.2.11)$$
on d\'eduit la suite exacte
longue 2.2.9 en homologie :
$\delta^*_{N,M}$ est le
connectant de 2.2.11.
\CQFD\medskip

L'int\'er\^{e}t essentiel des suites
exactes longues 1.2.2
et 2.2.9 par rapport \`a
1.1.2 est que sous des hypoth\`eses tr\`es
naturelles,
elles ont la vertue
de se scinder en suites exactes courtes
gradu\'ees.
\medskip

\goodbreak
{\bf\S2.3 Preuve du Th\'eor\`eme 1
({\it cf.} 1.2.3).} 

Donnons d'abord des crit\`eres
(n\'ecessaires \`a la preuve du
Th\'eor\`eme 1)
pour casser les suites exactes
longues 1.2.2 et 2.2.9 en suites
exactes courtes gradu\'ees.
\smallskip

\proclaim 2.3.1 Proposition 3.
Si $M$ est (isomorphe \`a) un
$A\otimes B^o$-facteur direct
de $N$, nous avons un
isomorphisme de modules gradu\'es :
$$HH^*\big[{}^A_{\,0}
{}^{\,M\oplus N}_{\quad B}\big]
\cong HH^*\big[{}^A_{\,0}{}^{N}_{B}\big]
\times Ext^{*-1}_{\!A\otimes B^o}(N,M)
\times Ext^{*-1}_{\!A\otimes B^o}(M,M)
\times Ext^{*-1}_{\!A\otimes B^o}(M,N)\ .
\leqno\!\!(2.3.2)$$

\noindent{\bf 2.3.3 Preuve.}
Supposons que l'on ait $M=M'\oplus N$,
pour des $A\otimes B^o$-modules
$M'$ et $N$ : nous avons une
d\'ecomposition en somme directe
de complexes
$$C^*\{M'\oplus N,N\}=\Sigma^*\{M',N\}
\oplus {\cal C}^*\{M'\!,N,N\}\leqno\!\!(2.3.4)$$
(m\^eme d\'emonstration que celle pour
la d\'ecomposition 2.2.6)
o\`u ${\cal C}^*\{M'\!,N,N\}$ est le c\^one
du morphisme de complexes :
$$\matrix{
{}^{\quad C^*(A)}_{C^*_{tri}(M'\!,M')}\cr
{}^{C^*_{tri}(N,N)}_{C^*_{tri}(N,N)}\cr
{}^{C^*(B)}\cr}
\buildrel{\displaystyle \lambda^*_{M'\!,N,N}}
\over\longlongrightarrow \matrix{
{}^{C^*(A,End_{B^o}M')}_{\,\,C^*(A,End_{B^o}N)}\cr
\,{}^{C^*(A,End_{B^o}N)}_{C^*(B,End_{A}M')}\cr
{}^{C^*(B,End_{A}N)}_{C^*(B,End_{A}N)}\cr
}\quad,\leqno\!\!(2.3.5)$$
dont les composantes \'eventuellement
non nulles sont les $12$ morphismes naturels :
$\alpha_{M'}^*$, $2$ fois $\alpha_{N}^*$
(il y a $2$ fois $C^*(A,End_{B^o}N)$),
$i_{M'}^*=(i_{M'}^{A*},i_{M'}^{B^o*})$,
$2$ fois $i_{N}^*=(i_{N}^{A*},i_{N}^{B^o*})$ et,
enfin, $\beta_{M'}^*$ et $2$ fois $\beta_{N}^*$.
La remarque que nous faisons est
que la restriction
$\rho^*_{M',N,N}$ du morphisme
$\rho^*_{N,M'\oplus N}$
\`a ${\cal C}^*\{M',N,N\}$ est encore
surjectif (\'evident) et,
plus encore, il existe un
morphisme de complexes
$\sigma^*_{M',N,N}:
{\cal C}^*\{M',N\}\longrightarrow
{\cal C}^*\{M',N,N\}$ tel que
$\rho^*_{M',N,N}\circ\sigma^*_{M',N,N}=1$,
{\it i.e.} une section, d\'efinie
de mani\`ere \'evidente
({\it cf.} 2.2.7 en rempla\c cant
$M$ par $M'$) en envoyant
sur eux-m\^emes
$C^*(A)$,
$C^*(B)$, $C^*_{tri}(M',M')$,
$C^*(A,End_{B^o}M')$ et
$C^*(B,End_{A}M')$ et en
envoyant $C^*_{tri}(N,N)$
(ainsi que $C^*(A,End_{B^o}N)$
et $C^*(B,End_{A}N)$)
sur ses deux copies
\`a l'aide
du morphisme diagonal.
\CQFD\medskip

La Proposition 3 exprime l'exactitude
et le scindage de la suite gradu\'ee
({\it cf.} 1.2.2) :
$$0\longrightarrow
\matrix{Ext^{*-1}_{A\otimes B^o}(N,M)\cr
Ext^{*-1}_{A\otimes B^o}(M,M)\cr
Ext^{*-1}_{A\otimes B^o}(M,N)\cr}
\buildrel{l^*_{N,M}}\over\longlongrightarrow
HH^*[{}^A_{\,0}\,{}^{M\oplus N}_{\ \ B}]
\buildrel{r^*_{N,M}}\over\longlongrightarrow
HH^*[{}^A_{\,0}{}^N_B]\longrightarrow 0\ ,
\leqno\!\!(2.3.6)$$
impossible \`a obtenir
avec 1.1.2 ($N\!=\!0$), \`a moins
d'avoir la situation triviale
($M\!=\!0$).

L'\'enonc\'e suivant est un renforcement
de la Proposition 3.

\proclaim 2.3.7 Proposition 4. Si $M$ est
(isomorphe \`a) un $A\otimes B^o$-facteur
direct de $N^m$, pour un certain $m\in{\bf N}$,
la suite gradu\'ee 2.3.6 est exacte et scind\'ee.

Une autre fa\c con d'\'enoncer ce r\'esultat
est la suivante. Posons
$A'=End_{A\otimes B^o}N$, ainsi
$N$ est un $A\otimes B^o\otimes A'$-module.
L'hypoth\`ese sur $M$ implique que
$Hom_{A\otimes B^o}(N,M)$ est un
$(A')^o$-module
projectif de type fini et, donc,
que $M'=Hom_{A\otimes B^o}(N,M)\oplus A'$
est un $(A')^o$-prog\'en\'erateur.
R\'eciproquement, 
tout $(A')^o$-prog\'en\'erateur de
la forme $M'=M''\oplus A'$,
v\'erifie $M''=Hom_{A\otimes B^o}(N,M)$
pour un certain facteur direct $M$
de $N^n$
(en fait, $M=M''\otimes_{A'}N$,
{\it cf. \bf[La]}, p.500).
De l\`a, la reformulation
suivante de la Proposition 4,
plus conceptuelle.

\proclaim 2.3.8 Proposition 4 (bis).
Pour tout $(A')^o$-prog\'en\'erateur
de la forme $M'=M''\oplus A'$,
nous avons un isomorphisme gradu\'e :
$$\eqalign{
H\!H^*\Big[{}^{A\ M'\otimes_{A'}N}_{\,0\quad\ B}\Big]
\cong &\quad H\!H^*[{}^{A\, M}_{\,0\, B}]\times
Ext^{*-1}_{A\otimes B^o}(N,N)\cr
&\times Ext^{*-1}_{A\otimes B^o}
(N,M''\!\otimes_{\!A'}\!N)
\times Ext^{*-1}_{A\otimes B^o}
(M''\!\otimes_{\!A'}\!N,N)\cr}
\ .\leqno\!\!(2.3.9)$$

\noindent{\bf 2.3.10 Preuve de la Proposition 4.}
Nous avons un diagramme commutatif :
$$\xymatrix{
HH^*[{}^A_{\,0}\,{}^{M\oplus N^m}_{\quad B}]
\ar[rr]^{r^*_{M\oplus N,N^{m-1}}}
\ar[d]_{r^*_{N^m,M}}
&&HH^*[{}^A_{\,0}\,{}^{M\oplus N}_{\quad B}]
\ar[d]^{r^*_{N,M}\qquad\textstyle .}\cr
HH^*[{}^A_{\,0}\,{}^{N^m}_{\,B}]
\ar[rr]^{r^*_{N,N^{m-1}}}
&&HH^*[{}^A_{\,0}\,{}^{N}_{B}]\cr}
\leqno\!\!(2.3.11)$$
D'apr\`es 2.3.6, les morphismes
$r^*_{N^m,M}$, $r^*_{M\oplus N,N^{m-1}}$
et $r^*_{N,N^{m-1}}$ sont surjectifs
et admettent des sections. Nous en
d\'eduisons la m\^eme chose
pour le morphisme $r^*_{N,M}$.
\CQFD\medskip

\noindent{\bf 2.3.12 Preuve du Th\'eor\`eme 1.}
Avec les notations du Th\'eor\`eme 1,
posons $M=\tilde M\oplus \tilde N$, avec
$\tilde M=\oplus_i M_i^{m_i-1}$
et $\tilde N=\oplus_i M_i$ (possible
si $m_i\geq 1$) ; ainsi $\tilde M$
est un facteur direct de $\tilde N^{\tilde m}$
pour $\tilde m=\max_i\{m_i-1\}$.
Ainsi, la Proposition 4 assure
l'exactitude et le scindage de 2.3.6
avec $\tilde M$ et $\tilde N$
\`a la place de $M$ et $N$,
c'est-\`a-dire l'exactitude
et le scindage de 1.2.4.
\CQFD\medskip

Nous terminons ce paragraphe avec
une ``formule d'\'echange'',
cons\'equence directe de la
Proposition 4.

\proclaim 2.3.13 Proposition 5. Pour
deux $A\otimes B^o$-modules $M$ et $N$
($K$-projectifs), tels que $M$
soit un facteur direct d'un $N^m$ et $N$
soit un facteur direct
d'un $M^n$ ($m,n\in{\bf N}$),
nous avons :
$$HH^*[{}^A_{\,0}\,{}^{M}_{\,B}]
\times Ext^{*-1}_{A\otimes B^o}(N,N)
\cong HH^*[{}^A_{\,0}\,{}^{N}_{B}]
\times Ext^{*-1}_{A\otimes B^o}(M,M)
\ .\leqno\!\!(2.3.14)$$
\bigskip

\goodbreak
\noindent\hfil{\bf3.
Second Th\'eor\`eme de R\'eduction}
\smallskip

Dans cette section, nous montrons
la suite exacte de Mayer-Vietoris
(le Th\'eor\`eme 2 de {\bf\S1.2}).
Dans {\bf\S3.1}, nous explicitons
le complexe d\'efinissant la cohomologie
modifi\'ee ${\cal H}^*$ et
nous montrons le Lemme 2 et
le Lemme 3 de {\bf\S1.2}.
La suite exacte
de mayer-Vietoris est obtenue
dans {\bf\S3.2}, comme cas particulier
d'un r\'esultat plus g\'en\'eral
(Th\'eor\`eme 4),
dont nous tirons d'autres
conclusions int\'eressantes,
telles les suites exactes 3.2.3 et 3.2.5.
\smallskip

\goodbreak
{\bf\S3.1 Cohomologie modifi\'ee.}

Soient deux alg\`ebres $A$ et $B$
($K$-modules projectifs) et soit
une famille de $A\otimes B^o$-modules
($K$-projectifs) $E=\{M_1,M_2,...,M_n\}$
($n\geq 0$) :
nous leur associons une alg\`ebre
triangulaire $T_E=[{}^A_{\,0}{}^{\bar E}_B]$
($\bar E=\oplus_iM_i$).
Avec les notations de {\bf\S2.1},
introduisons le complexe
$$\Sigma^*E=\prod_{i\not=j}
\Sigma^*_{M_i,M_j}\ .\leqno\!\!(3.1.1)$$
($\Sigma^*\phi=0$, par convention).
Le module de cohomologie de $\Sigma^*E$
est donn\'e en 2.1.4 :
$$H(\Sigma^*E)\cong \prod_{i\not=j}
Ext^{*-1}_{\!A\otimes B^o}(M_i,M_j)\ .
\leqno\!\!(3.1.2)$$
Aussi important est le morphisme
de complexes :
$$C^*(A)\oplus\bigoplus_iC^*_{tri}(M_i)
\oplus C^*(B)
\buildrel{\displaystyle\lambda^*_E\ }
\over\longrightarrow
 \bigoplus_i C^*(A,End_{B^o}M_i)
\oplus  \bigoplus_i C^*(B,End_{A}M_i)
,\leqno\!\!(3.1.3)$$
dont les seules composantes
\'eventuellement non nulles
sont repr\'esent\'ees ci-dessous
$$\xymatrix@R=2pt{
{\scriptstyle C^*(A)}\ar[rrr]\ar[rrrd]
\ar@{-->}[rrrdd]\ar[rrrddd]&&&\cr
{\scriptstyle C^*_{tri}(M_1)}\ar[rrru]
\ar[rrrddd]&&&\cr 
{\scriptstyle C^*_{tri}(M_2)}\ar[rrru]
\ar@{-->}[rrrddd]&&&\cr 
:&&&\cr 
:&&&\cr 
{\scriptstyle C^*_{\!tri}(\!M_{n-1}\!)}
\ar[rrrd]
\ar@{-->}[rrruuu]&&&\cr 
{\scriptstyle C^*_{tri}(M_n)}\ar[rrrd]
\ar[rrruuu]&&&\cr 
{\scriptstyle C^*(B)}\ar[rrr]
\ar[rrru]\ar@{-->}[rrruu]
\ar[rrruuu]&&&\cr}
\!\!\!\!\xymatrix@R=2pt{
{\scriptstyle C^*(A,End_{B^o}M_1)}\cr
{\scriptstyle C^*(A,End_{B^o}M_2)}\cr 
\qquad :\qquad\cr 
{\scriptstyle C^*(A,End_{B^o}M_n)}\cr 
{\scriptstyle C^*(B,End_{A}M_1)}\cr 
\qquad:\qquad\cr 
{\scriptstyle C^*(B,End_{A}M_{n-1})}\cr 
{\scriptstyle C^*(B,End_{A}M_n)}\cr}
\leqno\!\!(3.1.4)$$
Ce qui constitue un ensemble
de $4n$ fl\`eches : les $n$ morphismes
$\alpha_{M_i}^*$, les $n$ morphismes
$\beta_{M_i}^*$ et les $n$ morphismes
$i_{M_i}^{*}=(i_{M_i}^{A*},i_{M_i}^{B^o\!*})$
($2$ fl\`eches par morphisme).
Le morphisme $\lambda^*_E$
g\'en\'eralise la D\'efinition 3
($E=\{M\}$ dans 2.1.6) et l'on retrouve
2.2.7 ({\it resp.} 2.3.5) comme cas
particulier : $E=\{M,N\}$
({\it resp.} $E=\{M',N,N\}$).
La d\'efinition
suivante est donc attendue.
\proclaim 3.1.5 D\'efinition 4.
Nous posons
${\cal C}^*E=c\hat one(\lambda_E^*)$ et
la cohomologie modifi\'ee de $T_E$,
not\'ee ${\cal H}^*E$,
est la cohomologie de cohomologie de
${\cal C}^*E$.

\noindent{\bf 3.1.6 Preuve du Lemme 2
({\it cf.} 1.2.5).} 
Montrons que nous avons un
isomorphisme de complexes
${\cal C}^*\{\!\bar E\}\!=\!
{\cal C}^*E\!\oplus\! \Sigma^*E$ ;
nous aurons ainsi, d'apr\`es 3.1.1,
l'isomorphisme gradu\'e du Lemme 2 :
$$HH^*(T_E)\cong {\cal H}^*E\oplus
\prod_{i\not=j}
Ext^{*-1}_{\!A\otimes B^o}(M_i,M_j)
\ ,\leqno\!\!(3.1.7)$$
ce qui assure la compatibilit\'e
des deux d\'efinitions, celle ci-dessus
et la D\'efinition 1 donn\'ee dans {\bf\S1.2}.
D'apr\`es 2.1.5, le complexe
${\cal C}^*\{\bar E\}$ est le c\^one
du morphisme $\lambda^*_{\bar E}$
et nous avons des isomorphismes
(d\'ecompositions de Peirce avec $n$ idempotents) :
$$End_{B^o}\bar E=\prod_{i,j}Hom_{B^o}(M_i,M_j)
\quad{\rm et}\quad
End_A\bar E=\prod_{i,j}Hom_A(M_i,M_j)\ .$$
Dans ces d\'ecompositions, les
morphismes canoniques d'alg\`ebres
$\alpha_{\bar E}$ et $\beta_{\bar E}$ 
v\'erifient
$$Im\, \alpha_{\bar E}\subset
\prod_{i}Hom_{B^o}(M_i,M_i)
\quad{\rm et}\quad
Im\, \beta_{\bar E}\subset
\prod_{i}Hom_{A}(M_i,M_i)\ .$$
Nous en d\'eduisons la m\^eme
chose aux niveaux des complexes
de Hochschild :
$$Im\, \alpha^*_{\bar E}\subset
\prod_{i}C^*(A,End_{B^o}M_i)
\quad{\rm et}\quad
Im\, \beta^*_{\bar E}\subset
\prod_{i}C^*(B,End_{A}M_i)\ .$$
Le morphisme $\alpha^*_{\bar E}$
({\it resp.} $\beta^*_{\bar E}$) se r\'eduit
\`a ses composantes diagonales
$\alpha^*_{M_i}$
({\it resp.} $\beta^*_{M_i}$).
Nous en d\'eduisons que le
sous-complexe de ${\cal C}^*\{\bar E\}$
form\'e des $C^*_{tri}(M_i,M_j)$,
$C^*(A,Hom_{B^o}(M_i,M_j))$ et
$C^*(B,Hom_{A}(M_i,M_j))$ ($i\!\not=\!j$),
qui est naturellement isomorphe \`a $\Sigma^*E$,
est un facteur direct
de ${\cal C}^*\{\!\bar E\}$ et admet pour
complexe suppl\'ementaire le
complexe ${\cal C}^*E$
d\'efini en 3.1.5.
\CQFD\medskip

\noindent
{\bf 3.1.9 Preuve du Lemme 3 ({\it cf.} 1.2.7).}
Avec les notations pr\'ec\'edentes,
consid\'erons une partie $F\subset E$.
Il lui est associ\'e un morphisme
de complexes $\lambda^*_F$ (3.1.3),
un complexe ${\cal C}^*F$ et une
cohomologie ${\cal H}^*F$ (3.1.5).
D'apr\`es 3.1.4, le complexe ${\cal C}^*F$
est naturellement un quotient
de ${\cal C}^*E$
(rarement un sous-complexe).
Le noyau de la projection,
form\'e des sous-modules
$C^*_{tri}(L)$,
$C^*(A,End_{B^o}L)$
et $C^*(B,End_{A}L)$
avec $L\in E\backslash F$,
est un sous-complexe isomorphe
\`a $\prod_{L\in E\backslash F}\Sigma^*_L$.
Nous en d\'eduisons une suite exacte
courte de complexes :
$$0\longrightarrow
\prod_{L\in E\backslash F}\Sigma^*_L
\buildrel{\ell_{F,E}}
\over\longrightarrow {\cal C}^*E
\buildrel{\rho_{F,E}}
\over\longrightarrow{\cal C}^*F
\longrightarrow 0\ ,\leqno\!\!(3.1.10)$$
\`a l'origine de la suite exacte
longue 1.2.8 ; $\delta^*_{F,E}$
est le morphisme connectant de la suite
exacte courte 3.1.10.
Notant que cette derni\`ere est une
sous-suite exacte courte de 2.2.11
(o\`u l'on prend $N=\bar F$ et
$M=\overline{E\backslash F}$), il
est \`a pr\'esent clair que les
morphismes $\ell^*_{F,E}$, $\rho^*_{F,E}$
et $\delta^*_{F,E}$ sont 
induits, respectivement,
par $l^*_{\bar F,\bar E}$,
$r^*_{\bar F,\bar E}$ 
et $d^*_{\bar F,\bar E}$
de 1.2.2.
\CQFD\medskip

\noindent{\bf 3.1.11 Remarque.}
Pour $F\subset E$, notons que la
projection $\rho_{F,E}$ poss\`ede
une section canonique,
{\it a priori} seulement lin\'eaire,
$\sigma_{E,F}:
{\cal C}^*F\longrightarrow{\cal C}^*E$.
Elle aussi est fonctorielle :
$\sigma_{E,F}\circ \sigma_{F,G}=\sigma_{E,G}$,
pour $G\subset F\subset E$,
et, pour $F,G\subset E$ quelconques,
nous avons :
$$\rho_{F,E}\circ \sigma_{E,G}=
\sigma_{F,F\cap G}\circ
\rho_{F\cap G,G}\ .\leqno\!\!(3.1.12)$$
\medskip

\goodbreak
{\bf\S3.2 Suite exacte
de Mayer-Vietoris.}

Soit ${\cal U}=\{U_0,U_1,...,U_u\}$ un
recouvrement de $E$ ({\it i.e.} $\bigcup_iU_i=E$).
\proclaim 3.2.1 Th\'eor\`eme 4.
Supposons $U_i\cap U_j\subset U_0$
si $i\not=j$ ; nous avons une
suite exacte longue
$$\cdots\longrightarrow
{\cal H}^*E
\buildrel{\rho^*_+}
\over\longrightarrow
\prod_{i=0}^u{\cal H}^*U_i
\buildrel{\rho^*_-}
\over\longrightarrow
\prod_{i=1}^u{\cal H}^*(U_{0}\cap U_{i})
\longrightarrow{\cal H}^{*+1}E
\longrightarrow\cdots\leqno\!\!(3.2.2)$$
o\`u
$\rho^*_+=(\rho^*_{U_i,E}
)_{0\leq i\leq u}$
et $\rho^*_-=(
\rho^*_{U_{0}\cap U_{i},U_0}-
\rho^*_{U_{0}\cap U_{i},U_{i}}
)_{1\leq i\leq u}$.

\noindent \`A titre d'exemple
($U_0\!=\!\{M_1\}$ et
$U_i\!=\!\{M_1,M_{i+1}\}$ si $i\!>\!1$),
nous avons la suite
exacte longue :
$$\cdots\longrightarrow
\matrix{\prod_{i=2}^n
Ext^{*-1}_{\!A\otimes B^o}(M_1,\!M_i)\cr
{\cal H}^*E\cr
\prod_{i=2}^n
Ext^{*-1}_{\!A\otimes B^o}(M_i,\!M_1)\cr}
\longrightarrow
\prod_{i=2}^nHH^*[{}^A_{\,0}
{}^{M\!{\scriptscriptstyle 1}\oplus
M\!\scriptscriptstyle i}_{\quad B}]
\longrightarrow
(HH^*[{}^A_{\,0}
{}^{M\!\scriptscriptstyle 1}_{\ B}])^{n-2}
\longrightarrow\cdots\leqno\!\!(3.2.3)$$
D'autres cas particuliers m\'eritent
notre attention : si $u=1$ 
on retrouve le Th\'eor\`eme 2
(Mayer-Vietoris)
et, d'autre part, le cas suivant est important
en pratique.
\proclaim 3.2.4 Corollaire 2.
Si ${\cal U}$ est une partition de $E$,
nous avons une suite exacte longue :
$$\cdots\longrightarrow
{\cal H}^*E
\longrightarrow
\prod_{i=0}^u{\cal H}^*U_i
\longrightarrow
\big(HH^*A\times HH^*B\big)^u
\longrightarrow{\cal H}^{*+1}E
\longrightarrow\cdots\leqno\!\!(3.2.5)$$

\noindent{\bf 3.2.6 Preuve du Th\'eor\`eme 4.}
Consid\'erons le complexe
$\prod_{i=0}^u{\cal C}^*U_i$ : il est
muni d'un morphisme canonique
$\rho_+:{\cal C}^*E\longrightarrow
\prod_{i=0}^u{\cal C}^*U_i$ avec,
pour composantes, les
morphismes $\rho_{U_i,E}$
($0\leq i\leq u$) d\'efinis par 3.1.10.
Le fait que ${\cal U}$ soit un recouvrement
assure l'injectivit\'e de $\rho_+$.
Le morphisme $\rho_-$ est fabriqu\'e \`a partir
des morphismes surjectifs
$\rho_{U_0\cap U_i,U_{0}}$ et
$-\rho_{U_0\cap U_i,U_i}$ ($i>0$)
selon le sch\'ema d'assemblage
suivant :
$$\xymatrix@C=3mm@R=6mm{
{\cal C}^*U_0\ar[rd]\ar[rrd]
\ar@{-->}[rrrd]\ar[rrrrd]&
{\cal C}^*U_1\ar[d]&
{\cal C}^*U_2\ar[d]&
\quad.\,.\,.\,.\,.\quad\ar@{-->}[d]&
{\cal C}^*U_u\ar[d]\cr
&{\cal C}^*(U_0\cap U_1)
&{\cal C}^*(U_0\cap U_2)
&\quad.\,.\,.\,.\,.\quad
&{\cal C}^*(U_0\cap U_u)\cr}\leqno\!\!(3.2.7)$$
Ce qui prouve que $\rho_-$ est surjectif.
Montrons que nous avons
$Ker \rho_-=Im \rho_+$ : nous en d\'eduirons
une suite exacte courte de complexes
\`a l'origine de la suite exacte longue 3.2.5.
L'inclusion $Im\,\rho_+\subset Ker\,\rho_-$
\'etant claire avec
les formules de fonctorialit\'e
($\rho_{U_0\cap U_i,U_i}\circ \rho_{U_i,E}
=\rho_{U_0\cap U_i,E}$),
montrons l'autre.
Le morphisme
$\rho_-$ s'annule sur
$\prod_iC_i\subset \prod_i{\cal C}^*U_i$,
avec $C_i=Ker\,\rho_{U_0\cap U_i,U_i}$
si $i>0$ et
$C_0=Ker\,\rho_{U'_0,U_0}$,
avec $U'_0=U_0\cap \bigcup_{j>0}U_j$
(pour ce dernier :
$\rho_-(y_0,0,..,0)=\rho'_+(\rho_{U'_0,U_0}(y_0))$,
o\`u $\rho'_+$ est l'injection 
associ\'ee au recouvrement
$\{U_0\cap U_1,...,U_0\cap U_u\}$
de $U'_0$). Nous avons aussi
$Ker\,\rho_-=\prod_iC_i\oplus Ker\,\rho'_0$,
o\`u $\rho'_-$ est induit par $\rho_-$
({\it cf.} 3.1.11):
$$Im\,\sigma_{U_0,U'_0}
\times \prod_{i>0} Im\,\sigma_{U_i,U_0\cap U_i}
\longrightarrow
\prod_i{\cal C}^*(U_0\cap U_i)\ .
\leqno\!\!(3.2.8)$$
D'une part, les images des
$C_i$ par les sections
$\sigma_{U_i,E}$ sont disjointes
puisque
$U_i\cap U_j\subset U_0$
pour $i\not=j$ et
nous en d\'eduisons
$\prod_iC_i\subset Im\,\rho_+$.
D'autre part la restriction
de $\rho'_-$ \`a $\{0\}\times
\prod_{i>0} Im\,\sigma_{U_i,U_0\cap U_i}$
\'etant un isomorphisme, tout
\'el\'ement de $Ker\,\rho'_-$
s'\'ecrit :
$$y=(\sigma_{U_0,U'_0}x,
\sigma_{U_1,U_0\cap U_1}\rho_{U_0\cap U_1,U'_0}x,
\sigma_{U_2,U_0\cap U_2}\rho_{U_0\cap U_2,U'_0}x,
...,\sigma_{U_u,U_0\cap U_u}
\rho_{U_0\cap U_u,U'_0}x)\ ,$$
avec $x\in {\cal C}^*U'_0$.
Posant $z=\sigma_{E,U'_0}x$,
nous obtenons $y=\rho^+(z)$
en appliquant la formule 3.1.12
\`a chacune des coordonn\'ees
de $y$.
\CQFD\medskip

\goodbreak
\noindent\hfil{\bf4. L'Alg\`ebre de Lie HH$^1$T}
\smallskip

Pour deux alg\`ebres $A$, $B$ et
un bimodule $M$,
non suppos\'es $K$-modules projectifs,
nous \'etudions dans cette section
l'alg\`ebre de Lie
$HH^1\big[{}^{A\,M}_{\,0\ B}\big]$ et
son comportement relativement \`a
une d\'ecomposition de $M$.
Nous commen\c cons par
d\'emontrer le Th\'eor\`eme 3
de {\bf\S1.2} et, dans  {\bf\S4.2},
nous donnons des conditions suffisantes
pour l'obtention de la suite exacte
1.2.16.
\smallskip

{\bf\S4.1 Preuve du Th\'eor\`eme 3
({\it cf.} 1.2.13).}

Pour toute $K$-alg\`ebre $T$,
nous savons qu'il y a un isomorphisme
d'alg\`ebres de Lie $HH^1T\cong Der(T)/Int(T)$,
o\`u $Der(T)$ est l'alg\`ebre de Lie
des $K$-d\'erivations $D$ de $T$ :
\smallskip

$\forall\,k\in K,\ \forall\,t,t'\in T,\
D(tt')=D(t)t'+tD(t')$
et $D(kt)=kD(t)$ 
\smallskip

\noindent
et $Int(T)$ est l'id\'eal de Lie
des ``d\'erivations int\'erieures''
$[t_0,.]$, $t_0\in T$
($\forall\,t,\ [t_0,t]=t_0t-tt_0$).

Dans le cas de
l'alg\`ebre triangulaire,
nous avons :

\proclaim 4.1.1 Lemme 4. Tout \'el\'ement
$D\in Der[{}^{A\,M}_{\,0\ B}]$ est d\'ecrit
par $\alpha\in Der(A)$, $\beta\in Der(B)$,
$m_0\in M$ et $\mu\in End_K(M)$ :
$\forall [{}^{\,a\,m}_{\,0\ b}]\in
[{}^{A\,M}_{\,0\ B}]$,
$D[{}^{\,a\,m}_{\,0\ b}]=[{}^{\alpha(a)\
\mu(m)-am_0+m_0b\,
}_{\ 0\qquad\qquad \beta(b)}]$.
Ils v\'erifient :
$$\mu(am)=\alpha(a)m+a\mu(m)
\quad{\sl et }\quad
\mu(mb)=\mu(m)b+m\beta(b)\ .
\leqno\!\!(4.1.2)$$

Nous noterons
$D=[{}^{\alpha\,(\mu,m_0)
}_{\,\,0\quad\beta}]$.
\medskip

\noindent{\bf 4.1.3 Preuve.}
Il s'agit d'un simple (et long) calcul.
Le traitement complet de l'\'equation
$D(tt')=D(t)t'+tD(t')$ implique la r\'esolution
d'un syst\`eme de $9$
morphismes inconnus
$D_{ij}:T_i\longrightarrow T_j$,
o\`u $T_1=A$, $T_2=M$, $T_3=B$ :
voir, par exemple, {\bf[FM]}
(et aussi {\bf[Ch]}).
\CQFD\medskip

Pour $[{}^{\,a_0\,m_0}_{\,\,0\ \ b_0}]\in
[{}^{A\,M}_{\,0\ B}]$, nous avons
une d\'erivation (int\'erieure)
$D=[{}^{\alpha\,(\mu,m_0)
}_{\,\,0\quad\beta}]$ :
$$D[{}^{\,a\,m}_{\,0\ b}]=
\big[{}^{[a_0,a]\ a_0m-am_0+m_0b-mb_0
}_{\quad 0\ \qquad\qquad [b_0,b]}\big]
\ ,\leqno\!\!(4.1.4)$$
{\it i.e.} : 
$\alpha=[a_0,.]$ et $\beta=[b_0,.]$
(d\'erivations int\'erieures)
et $\mu(m)=a_0m-mb_0$, $m\in M$.
Nous en d\'eduisons une d\'ecomposition
en somme directe :
$$Int[{}^{A\,M}_{\,0\ B}]
\cong\Big((ZA\times ZB)/
\{(a,b)|\forall m\in M,\, am=mb\}\Big)
\oplus M\ ,\leqno\!\!(4.1.5)$$
o\`u $M$ est isomorphe
au sous-module $\{{}^{0\,M}_{0\ 0}\}$
des d\'erivations
de la forme
$\big[{}^{\,0\,(0,m_0)}_{\,0\quad 0}\big]$
($m_0\in M$).
Le lemme suivant est le r\'esultat
d'un calcul direct.
\proclaim 4.1.6 Lemme 5.
Pour deux d\'erivations
$D_0=[{}^{\alpha_0\,(\mu_0,m_0)
}_{\,0\quad\ \beta_0}],
D_1=[{}^{\alpha_1\,(\mu_1,m_1)
}_{\,0\quad\ \beta_1}]$,
nous avons :
$$[D_0,D_1]=\big[{}^{\,[\alpha_0,\alpha_1]\
([\mu_0,\mu_1],\mu_0(m_1)-\mu_1(m_0))
}_{\quad\ 0\qquad\qquad\quad
[\beta_0,\beta_1]\,}\big]
\ .\leqno\!\!(4.1.7)$$

\noindent En particulier,
$\{{}^{0\,M}_{0\ 0}\}$
est un id\'eal de Lie
(ab\'elien)
de $Der[{}^{A\,M}_{\,0\ B}]$ :
$Int'[{}^{A\,M}_{\,0\ B}]=Int[{}^{A\,M}_{\,0\ B}]
/\{{}^{0\,M}_{0\ 0}\}$
est un id\'eal de l'alg\`ebre
de Lie $Der'[{}^{A\,M}_{\,0\ B}]
=Der[{}^{A\,M}_{\,0\ B}]/\{{}^{0\,M}_{0\ 0}\}$
et nous avons un isomorphisme
d'alg\`ebres de Lie : 
$HH^1[{}^{A\,M}_{\,0\ B}]=
Der'[{}^{A\,M}_{\,0\ B}]/Int'[{}^{A\,M}_{\,0\ B}]$.
La sous-alg\`ebre de Lie 
de $Der[{}^{A\,M}_{\,0\ B}]$
form\'ee des \'el\'ements
$D=[{}^{\alpha\,(\mu,0)
}_{\,0\ \ \beta}]$
est isomorphe \`a
$Der'[{}^{A\,M}_{\,0\ B}]$ :
on \'ecrira $D=
[{}^{\alpha\,\mu}_{0\,\beta}]$.
\smallskip

Soient \`a pr\'esent
un second $A\otimes B^o$-module
$N$ et une
d\'erivation
$[{}^{\,\alpha\,\gamma
}_{\,0\,\,\beta}]\in 
Der'[{}^{A\,M\oplus N}_{\,0\quad B}]$.
\proclaim 4.1.8 Lemme 6.
Les composantes
$\mu\in End_K(M)$ et
$\nu\in End_K(N)$ de $\gamma$,
v\'erifient :
$$[{}^{\,\alpha\,\mu
}_{\,0\,\beta}]\in 
Der'[{}^{A\,M}_{\,0\ B}]
\quad{\rm et}\quad [{}^{\,\alpha\,\nu
}_{\,0\,\beta}]\in 
Der'[{}^{A\,M}_{\,0\ B}]
\ ;\leqno\!\!(4.1.9)$$
les composantes de $\gamma$
``non diagonales'' $M\rightarrow N$
et $N\rightarrow M$ sont
des $A\otimes B^o$-morphismes.

\noindent{\bf 4.1.10 Preuve.}
C'est un simple calcul.\CQFD\medskip

\noindent Nous
en d\'eduisons deux morphismes
(non n\'ecessairement
morphismes d'alg\`ebres de Lie) :
$$\eqalign{
r'_{M,N}&:Der'[{}^{A\,M\oplus N}_{\,0\quad B}]
\longrightarrow Der'[{}^{A\,M}_{\,0\ B}]\ ,\
[{}^{\,\alpha\,\gamma}_{\,0\ \beta}]
\longmapsto
[{}^{\,\alpha\,\mu}_{\,0\,\beta}]\ ,\cr
r'_{N,M}&:
Der'[{}^{A\,M\oplus N}_{\,0\quad B}]
\longrightarrow Der'[{}^{A\,N}_{\,0\, B}]\ ,\
[{}^{\,\alpha\,\gamma}_{\,0\ \beta}]
\longmapsto
[{}^{\,\alpha\,\nu}_{\,0\,\beta}]\ ,\cr}
\leqno\!\!(4.1.11)$$
v\'erifiant des relations de transitivit\'e
($r'_{N,M}r'_{M\oplus N,L\oplus M\oplus N}
\!=\!r'_{N,L\oplus M\oplus N}\!=\!
r'_{N,L}r'_{N\oplus L,L\oplus M\oplus N}$).
Puisque
$r'_{M,N}
\big(Int'[{}^{A\,M\oplus N}_{\,0\quad B}]\big)
\subset Int'[{}^{A\,M}_{\,0\ B}]$ et
$r'_{N,M}
\big(Int'[{}^{A\,M\oplus N}_{\,0\quad B}]\big)
\subset Int'[{}^{A\,N}_{\,0\ B}]$,
ils d\'efinissent
des morphismes (non morphismes d'alg\`ebres
de Lie {\it a priori}) :
$$r^1_{M,N}:HH^1[{}^{A\,M\oplus N}_{\,0\quad B}]
\longrightarrow HH^1[{}^{A\,M}_{\,0\ B}]
\quad{\rm et}\quad
r^1_{N,M}:HH^1[{}^{A\,M\oplus N}_{\,0\quad B}]
\longrightarrow HH^1[{}^{A\,N}_{\,0\ B}]
\ .\leqno\!\!(4.1.12)$$

Du Lemme 5, nous d\'eduisons une
d\'ecomposition en somme directe :
$$Der'[{}^{A\,M\oplus N}_{\,0\quad B}]=
D'\{M,N\}\oplus D^o(M,N)\oplus D^o(N,M)\ ,$$
o\`u  $D'\{M,N\}=D'\{N,M\}$ est la sous-alg\`ebre
de Lie des \'el\'ements
$[{}^{\,\alpha\,\gamma
}_{\,0\,\,\beta}]\in
Der'[{}^{A\,M\oplus N}_{\,0\quad B}]$,
tels que les composantes
non diagonales
de $\gamma$ sont nulles et
$D^o(M,N)\cong Hom_{A\otimes B^o}(M,N)$
(resp. $D^o(N,M)\cong Hom_{A\otimes B^o}(N,M)$)
est l'id\'eal (ab\'elien) form\'e
des \'el\'ements
$[{}^{\,0\,\gamma
}_{\,0\,0}]\in
Der'[{}^{A\,M\oplus N}_{\,0\quad B}]$,
o\`u $\gamma$ ne comporte au plus
qu'une composante non nulle
$M\longrightarrow N$
(resp. $N\longrightarrow M$).
Notons que l'alg\`ebre de Lie  $D'\{M,N\}$
est aussi l'intersection de deux
autres sous-alg\`ebres de Lie
(matrices ``triangulaires'' sup\'erieures
et inf\'erieures) :
$$D'(M,N)= D'\{M,N\}\oplus D^o(M,N)
\ {\rm et}\
D'(N,M)= D'\{M,N\}\oplus D^o(N,M)
\ .\leqno\!\!(4.1.13)$$

\noindent{\bf 4.1.14
Fin de la preuve du Th\'eor\`eme 3.}
Puisque $Int'[{}^{A\,M\oplus N}_{\,0\quad B}]
\subset  D'\{M,N\}$, nous en d\'eduisons
les trois sous-alg\`ebres de Lie de
$HH^1[{}^{A\,M\oplus N}_{\,0\quad B}]$ 
suivantes :
$$\eqalign{{\cal H}^1\{M,N\}&=
D'\{M,N\}/Int'
[{}^{A\,M\oplus N}_{\,0\quad B}]\ ,\cr
{\cal H}^1(M,N)&={\cal H}^1\{M,N\}\oplus 
D^o(M,N)\ ,\cr
{\cal H}^1(N,M)&={\cal H}^1\{M,N\}\oplus 
D^o(N,M)\ ,\cr}\leqno\!\!(4.1.15)$$
v\'erifiant 1.2.14.
Par restriction de $r'_{M,N}$,
nous obtenons plusieurs morphismes
$$\eqalign{
\rho'_{M,\{M,N\}}&:D'\{M,N\}\longrightarrow
Der'[{}^{A\,M}_{\,0\, B}]\ ,\cr
\rho'_{M,(M,N)}&:D'(M,N)\longrightarrow
Der'[{}^{A\,M}_{\,0\, B}]\ ,\cr
\rho'_{M,(N,M)}&:D'(N,M)\longrightarrow
Der'[{}^{A\,M}_{\,0\, B}]\ ,\cr}
\leqno \!\!(4.1.16)$$
ayant m\^eme image que  $r'_{M,N}$.
D'apr\`es 4.1.7 et 4.1.13,
il s'agit l\`a de morphismes
d'alg\`ebres de Lie, dont nous d\'eduisons
les morphismes d'alg\`ebres de Lie
1.2.15.
\CQFD\medskip

\goodbreak
{\bf\S4.2 Scindage de HH${}^1$T.}

Soient deux alg\`ebres $A$ et $B$
et soit un $A\otimes B^o$-module $M$
(nous ne supposons pas qu'ils sont
$K$-projectifs).
Nous notons $\delta[M]$ la
sous-cat\'egorie pleine des
$A\otimes B^o$-modules $N$,
tels que 
$r^1_{M,N}:HH^1[{}^{A\,M\oplus N}_{\,0\quad B}]
\longrightarrow HH^1[{}^{A\,M}_{\,0\, B}]$
(ou, de mani\`ere \'equivalente,
le morphisme d'alg\`ebres de Lie
$\rho^1_{M,\{M,N\}}:{\cal H}^1\{M,N\}
\longrightarrow HH^1[{}^{A\,M}_{\,0\, B}]$)
est surjectif ({\it cf.} 4.1.12 et 4.1.16).

Notons que si $HH^1[{}^{A\,M}_{\,0\, B}]=0$,
la cat\'egorie $\delta[M]$ est celle
de tous les $A\otimes B^o$-modules.

La transitivit\'e de $r^1$,
encod\'ee par la commutativit\'e du diagramme
$$\xymatrix{
HH^1[{}^{A\,M\oplus M'\oplus M''}_{
\,0\qquad B}]
\ar[r]\ar[d]\ar[rd]&
HH^1[{}^{A\,M\oplus M'}_{\,0\quad\ B}]
\ar[d]\cr
HH^1[{}^{A\,M\oplus M''}_{\,0\quad\ B}]
\ar[r]
&HH^1[{}^{A\,M}_{\,0\ B}]\cr}
\!\!\leqno(4.2.1)$$
(pour des $A\otimes B^o$-modules $M'$ et $M''$),
montre que nous avons :
$$M'\oplus M''\in \delta[M]
\Longrightarrow
M'\in \delta[M]\ {\rm et}\
M''\in \delta[M]
\ ,\leqno\!\!(4.2.2)$$
{\it i.e.} $\delta[M]$ est
stable par facteur direct.

\proclaim 4.2.3 Th\'eor\`eme 5.
$N\in \delta[M]$ si et
seulement si, pour tous
$\alpha\in Der(A)$ et
$\beta\in Der(B)$, nous avons :
s'il existe $\mu\in End_K(M)$
v\'erifiant
4.1.2, il existe $\nu\in End_K(N)$ 
v\'erifiant 4.1.2. Auquel cas,
nous avons une suite exacte
(dont 1.2.16 se d\'eduit) :
$$0\longrightarrow
Z\big[{}^{A\,M\oplus N}_{\,0\quad B}\big]
\buildrel{r^0_{M,N}}\over\longrightarrow
Z\big[{}^{A\,M}_{\,0\, B}\big]
\longrightarrow\!\!\!\!
\matrix{H\!om_{A\otimes B^o}(M,\!N)\cr
End_{A\otimes B^o}N\cr
H\!om_{A\otimes B^o}(N\!,M)\cr}
\!\!\!\!\longrightarrow
H\!H^1\big[{}^{A\,M\oplus N}_{\,0\quad B}\big]
\buildrel{r^1_{M,N}}\over\longrightarrow
H\!H^1\big[{}^{A\,M}_{\,0\, B}\big]
\longrightarrow 0.$$

\noindent Nous en d\'eduisons les
deux exemples suivants :
$$M\in \delta[M]\quad{\rm et}\quad
A\otimes B^o\in\delta[M]
\ .\leqno\!\!(4.2.4)$$
Le premier est (\`a pr\'esent) trivial et,
pour $N=A\otimes B^o$ et tout $\mu$,
$\nu=\alpha\otimes id+id\otimes\beta$
convient.
\medskip

\noindent{\bf 4.2.5 Preuve.}
D'apr\`es 4.1.5, nous avons
une surjection
$\tilde\rho_{M,\{M,N\}}\!:\!
Int'[{}^{A\,M\oplus N}_{\,0\quad B}]
\!\longrightarrow\! Int'[{}^{A\,M}_{\,0\ B}]$
(induite par $\rho'_{M,\{M,N\}}$,
{\it cf.} 4.1.16).
Nous en d\'eduisons un isomorphisme
(lemme du serpent) :
$$Coker\,\rho'_{M,\{M,N\}}\cong
Coker\,\rho^1_{M,\{M,N\}}\ .$$
En particulier : nous avons
$N\in\delta[M]$ si et seulement si
$\rho'_{M,\{M,N\}}$ est surjectif.
Le th\'eor\`eme se d\'eduit alors
de 4.1.15 (avec $N$ \`a la place de $M$)
 et de la suite exacte :
$$\xymatrix@R=1mm{0\ar[r]&
End_{A\otimes B^o}N\ar[r]&
D'\{M,N\}\ar[rr]^{\rho'_{M,\{M,N\}}}&&
Der'[{}^{A\,M}_{\,0\, B}]\cr
&\nu'\ar@{|->}[r]&[{}^{\,0\,\gamma'
}_{\,0\ 0}]&&\cr}\leqno\!\!(4.2.6)$$
o\`u $\gamma'\in End_K(M\oplus N)$ admet
$\nu'$ pour seule composante \'eventuellement
non nulle. 
\CQFD\medskip

\noindent
Les exemples 4.2.4 servent de 
briques de base pour les constructions
d'objets de $\delta[M]$ plus compliqu\'es,
d\'efinies par les corollaires suivants.

\proclaim 4.2.7 Corollaire 3. Pour tout
$A\otimes B^o$-module $N$, nous avons :
$$N\in \delta[M]\Longrightarrow
\big(
\delta[N]\subset \delta[M]\ {\sl et}\
\delta[M]=\delta[M\oplus N]\big)
\ .\leqno\!\!(4.2.8)$$

\noindent En particulier
$\delta[\,0\,]\subset \delta[M]$
et, pour tout $N\in \delta[\,0\,]$,
nous avons $\delta[N]=\delta[\,0\,]$.

\proclaim 4.2.9 Corollaire 4.
Pour toute famille
$N_\lambda\in \delta[M]$,
$\lambda\in\Lambda$,
nous avons :
$$\bigoplus_{\lambda\in\Lambda}
N_\lambda\in \delta[M]\quad{\sl et}\quad
\prod_{\lambda\in\Lambda}
N_\lambda \in \delta[M]
\ .\leqno\!\!(4.2.10)$$

\noindent Donc, pour tout cardinal
$m$ (m\^eme infini), nous avons
$M^m\in \delta[M]$ et
$M^{(m)}\in \delta[M]$ : c'est 
une g\'en\'eralisation de 1.3.2
(pour $*=1$). Nous avons aussi
$(A\otimes B^o)^{(m)}\in \delta[M]$ :
avec 4.2.2, ceci prouve que tout
$A\otimes B^o$-module projectif
est objet de $\delta[M]$.
Si nous disposions de l'hypoth\`ese
de projectivit\'e pour $A$ et $B$,
alors tout $A\otimes B^o$-module
($K$-projectif) $N$ tel que
$Ext^1_{A\otimes B^o}(N,N)=0$
(en particulier un $A\otimes B^o$-module injectif)
v\'erifierait $N\in \delta[\,0\,]$
d'apr\`es 1.1.2 et, en particulier,
$N\in \delta[M]$ ({\it cf.} aussi 2.2.9).

\proclaim 4.2.11 Corollaire 5.
Soit un morphisme d'alg\`ebres
$C\longrightarrow End_{A\otimes B^o}M$
({\it i.e.}, $M$ est muni d'une structure de
$A\otimes B^o\otimes C$-module).
Pour tout $C^o$-module $X$
et tout $C$-module $Y$ :
$$X\otimes_CM\in \delta[M]
\quad{\sl et}\quad
Hom_C(Y,M)\in \delta[M]
\ .\leqno\!\!(4.2.12)$$
De m\^eme, nous avons :
$Tor^C_i(X,M)\in \delta[M]$
et $Ext^i_C(Y,M)\in \delta[M]$,
$i>0$.

\bigskip

\noindent\hfil{\bf5.
R\'ef\'erences Bibliographiques}
\smallskip

\item{\bf [ARS]} M. Auslander, I. Reiten
\& S.O. Smal\o.
{\sl Representation Theory of Artin Algebras.}
Cambridge studies in advanced mathematics 36.
Cambridge university press, 1995.

\item{\bf [BG1]} B. Bendiffalah \&
D. Guin, {\sl Cohomologie de diagrammes
d'alg\`ebres triangulaires.}
Bolet\'\i n de la Academia Nacional
de Ciencias. C\'ordoba,  Argentina 65, pp.61--71, 2000.

\item{\bf [BG2]} B. Bendiffalah \&
D. Guin, {\sl Cohomologie de l'alg\`ebre
triangulaire et applications.}
Journal of Algebra 282, pp.513--537, 2004.

\item{\bf [CE]} H. Cartan \&
S. Eilenberg, {\sl
Homological Algebra.}
Princeton Mathematical Series,
Princeton university press, 1956.

\item{\bf [Ch]} W-S. Cheung,
{\sl Lie derivations of triangular algebras.}
Linear and Multilinear Algebra 51-3,
pp.299--310, 2003.

\item{\bf [Ci]} C. Cibils,
{\sl Tensor Hochshild homology
and cohomology.}
Lecture Notes in Applied
and Pure Mathematics, Vol. 210, pp.35--51, 2000.

\item {\bf [FM]} B.E. Forrest
\& L.W. Marcoux, {\sl Derivations
of triangular Banach algebras.}
Indiana Univ. Math. Jour. 45-2,
pp.441--462, 1996.

\item {\bf [GMS]} E.L. Green, E.N. Marcos
\& N. Snashall,
{\sl The Hochschild cohomology ring of
a one point extension.}
Communication in Algebra 31-1,
pp.357--379, 2003.

\item {\bf [GS]} E.L. Green \& \O. Solberg,
{\sl Hochschild
cohomology rings and triangular rings.}
Beijing Norm. Univ. Press,
Representations of algebras I \& II,
pp.192--200, 2002.

\item {\bf [GG]} J.A. Guccione \&
J.J. Guccione, {\sl Hochschild
cohomology of triangular matrix algebras.}
ArXiv: math.KT/0104068v2, 2001.

\item {\bf [GAS]} F. Guil-Asensio \& M. Saor\'\i n,
{\sl The automorphism group and the Picard group
of a monomial algebra.}
Communication in Algebra 27-2,
pp.857--887, 1999.

\item {\bf [Hap]} D. Happel,
{\sl Hochschild cohomology of
finite dimensional algebras.} Springer
Lecture Notes in Math. 1404, pp.108--126, 1989.

\item {\bf [Har]} M. Harada,
{\sl Hereditary semi-primary rings and
tri-angular matrix rings.} Nagoya Math. J. 27,
pp.463--484, 1966.

\item {\bf [K]} B. Keller,
{\sl Derived invariance of higher
structures on the Hochschild complex.}
Preprint, 2003.

\item{\bf [La]} T.Y. Lam,
{\sl Lectures on modules and rings.}
Graduate texts in mathematics,
Springer-Verlag New-York, 1999.

\item{\bf [Lo]} J-L. Loday,
{\sl Cyclic Homology,} Second Edition. Grundlehren
der mathematischen
Wissenschaften,
Springer-Verlag, Berlin, 1998.

\item{\bf [MP]} S. Michelena \& M. I. Platzeck,
{\sl Hochschild Cohomology of Triangular
Matrix Algebra.}
Journal of Algebra, Vol.233,
pp.502--525, 2000.

\item{\bf [M]} B. Mitchell,
{\sl Theory of categories.}
Pure and Applied Mathematics XVII.
Academic Press, New York-London, 1965.

\item{\bf [S]} C. Strametz,
{\sl  The Lie algebra structure of
the first Hochschild cohomology group
for monomial algebras.}
Comptes Rendus Acad. Sci. Paris 334-9,
pp.733--738, 2002.
\end